\documentclass[a4paper,12pt]{article}
\usepackage[top=2.5cm,bottom=2.5cm,left=2.5cm,right=2.5cm]{geometry}
\usepackage{amssymb}
\usepackage{graphicx}
\usepackage{amsmath,amsthm,amssymb,lineno}
\usepackage{latexsym}
\usepackage{graphicx,booktabs,multirow}
\usepackage{tikz}
\usetikzlibrary{decorations.pathreplacing}
\usetikzlibrary{intersections}
\usepackage{graphicx,booktabs,multirow}
\usepackage{appendix}

\newtheorem{theorem}{Theorem}
\newtheorem{lemma}[theorem]{Lemma}
\newtheorem{corollary}[theorem]{\rm\bfseries Corollary}

\begin{document}
	\title{The path sequence of a graph}
	\author{Yirong Cai, Hanyuan Deng\thanks{Corresponding author:
			hydeng@hunnu.edu.cn}\\
		{\footnotesize College of Mathematics and Statistics, Hunan Normal University, Changsha, Hunan 410081, P. R. China} \\
	}
	\date{}
	\maketitle
	
	\begin{abstract}
		Let $P(G)=(P_{0}(G),P_{1}(G),\cdots, P_{\rho}(G))$ be the path sequence of a graph $G$, where $P_{i}(G)$ is the number of paths with length $i$ and $\rho$ is the length of a longest path in $G$. In this paper, we first give the path sequences of some graphs and show that the number of paths with length $h$ in a starlike tree is completely determined by its branches of length not more than $h-2$. And then we consider whether the path sequence characterizes a graph from a different point of view and find that any two graphs in some graph families are isomorphic if and only if they have the same path sequence.\\
		\noindent
		{\bf Keywords}: Path sequence; Starlike tree; Generalized starlike tree; Kite graph; Lollipop graph.
		\end{abstract}
		\maketitle
	
	\makeatletter
	\renewcommand\@makefnmark%
	{\mbox{\textsuperscript{\normalfont\@thefnmark)}}}
	\makeatother
	
	\section{Introduction}\label{intro}

Many graphs can be determined by certain parameters, such as spectrum [1,3,9,13], Laplacian spectrum [7,14,15], walk sequence [2,6,11,12], and higher order indices [10,16]. Also, the study on the structure of starlike trees aroused the interest of many scholars. For example, Omidi and Tajbakhsh [8] proved that the starlike trees can be determined by their Laplacian spectrum. Bu and Zhou [4] proved that all starlike trees whose maximum degree exceed four are determined by their Q-spectra. Rada and Araujo [10] showed that starlike trees which have equal $h$-connectivity index for all $h\geq0$ are isomorphic.
	
Lexicographical ordering of graphs by spectral moments had been used earlier to order graphs, and the ordering can also obtain sufficient and necessary conditions for determining a graph in a class of graphs. Stevanovi\'{c} [12] considered the partial order relationship $\preceq$ for two graphs $G$ and $H$, $G\preceq H$, if $M_{k}(G)\leq M_{k}(H)$ for each $k\geq0$ and $M_{k}(G)$ denote the number of closed walks with length $k$ in $G$ (or $k$-th spectral moment). At the same time, the starlike trees $S_{n}$ with $n$ vertices is replaced by $S(a_{1},a_{2},\cdots, a_{l})$, where $a_{i}$ denotes the number of branches of $S_{n}$ with length $i$, where $1\leq i\leq l$ and $l$ is the length of a longest branch. It turns out that the partial order $\preceq$ of starlike trees with a fixed number of vertices coincides with the shortlex order of sorted sequence of their branch lengths. Let $M(S_{n})= (M_{0}(S_{n}),M_{1}(S_{n}),\cdots,M_{t}(S_{n}),\cdots)$ be the closed walk sequence of $S_{n}$. The result in [12] can further state that two starlike trees are isomorphic, if and only if their closed walk sequences are equal.
	
Although the number of closed walks of a graph can be calculated from some of the results in [5,6,11], the closed walk sequence of a graph is an infinite number sequence. Based on this, we investigate whether a graph can be determined by the path sequence, which has simpler structure than the closed walk sequence. Here, the path sequence of a graph $G$ is $P(G)=(P_{0}(G),P_{1}(G),\dots, P_{\rho}(G))$, where $P_{i}(G)$ is the number of paths in $G$ with length $i$ and $\rho$ is the length of a longest path in $G$. For convenience, we may assume that $P_i(G)=0$ for $i>\rho$. Specially, $P_{0}(G)=|V(G)|$, $P_{1}(G)=|E(G)|$. In more detail, if two graphs have the same path sequence, whether are they isomorphic? Or, can we use the path sequences to distinguish a graph?

Note that there are two cospectral connected graphs, their characteristic polynomial is $x^6-7x^4-4x^3+7x^2+4x-1$, these two graphs can not be distinguished by the spectrum and the (closed) walk sequence. However, they have different path sequences. Certainly, it is impossible to distinguish all graphs by the path sequence. For example, $K_{1,3}$ and $K_3\cup K_1$ have the same path sequence $(4,3,3)$.
	
	In this paper, we first give the number of the path with length $h$ in some graphs, including complete graphs, complete bipartite graphs, cycles, paths, stars, starlike graphs, kite graphs and lollipop graphs and so on. It is also found that the graph in some graph families can be determined by its path sequence. In other words, two graphs in these graph families with the same path sequence are isomorphic.

\section{The path sequence of some graphs}

	All graphs considered in this paper are simple and connected. It is hard to compute the path sequence as already the length of the longest path in a graph is NP-hard, formulas or methods for calculating path sequences of some special graph classes will become useful. With the definition of the path sequence, we will give in this section the path sequence of some special graphs and find these graphs can be determined by the path sequence easily. In order to get the path sequence of a graph, we first have to consider the number of paths with a given length in this graph. Such paths can be classified according to the degree sequence of their vertices.

 For a graph $G$, $\rho=\rho(G)$ denotes the length of a longest path in $G$ and $\mathcal{P}_{h}(G)$ denotes the set of all paths with length $h$ in $G$. Let $K_{n}$ be the complete graph with $n$ vertices. Since any path of length $h$ ($1\leq h\leq \rho(K_{n})=n-1$) in $K_{n}$ is a sequence of $h+1$ vertices without direction, we can easily get the following result.

\begin{theorem}\label{t-2}
The number of all paths with length $h$ in $K_{n}$ is
		$$	P_{h}(K_{n})=|\mathcal{P}_{h}(K_{n})| =\dfrac{1}{2}\prod_{i=0}^{h}(n-i), \;\text{where}\;1\leq h\leq n-1.$$
And the path sequence of $K_n$ is
$$P(K_{n})=(n,\dfrac{n(n-1)}{2},\dfrac{n(n-1)(n-2)}{2},\cdots,\dfrac{1}{2}\prod_{i=0}^{h}(n-i),\cdots,\dfrac{1}{2}\prod_{i=0}^{n-1}(n-i)).$$
	\end{theorem}
	
Note that a complete graph is determined by the number of vertices and $P_{0}(G)=|V(G)|$, we can get
	
\begin{corollary}\label{c-3}
		Two complete graphs are isomorphic if and only if they have the same path sequence.
	\end{corollary}
	
Let $K_{n_{1},n_{2}}$ be a complete bipartite graph with $n_{1}+n_{2}$ vertices and $n_1\leq n_2$. Then $\rho(K_{n_{1},n_{2}})=2n_1-1=2\min\{n_{1},n_{2}\}-1$ for $n_1=n_2$ and $\rho(K_{n_{1},n_{2}})=2n_1=2\min\{n_{1},n_{2}\}$ for $n_1\neq n_2$. We directly consider the number of paths with length $h$ in $K_{n_{1},n_{2}}$, denoted by $P_{h}(K_{n_{1},n_{2}})$.
	
	\begin{theorem}\label{l-4}
		Let $K_{n_{1},n_{2}}$ be a complete bipartite graph with $n_{1}+n_{2}$ vertices and $n_1\leq n_2$, $1\leq h\leq \rho(K_{n_{1},n_{2}})$. Then number of paths with length $h$ in $K_{n_{1},n_{2}}$ is
		$$P_{h}(K_{n_{1},n_{2}})=\left\{\begin{array}{l}\dfrac{1}{2}\prod\limits_{i=0}^{\frac{h}{2}}(n_{1}-i)(n_{2}-i)\times\left ( \dfrac{1}{n_{1}-\dfrac{h}{2}} +\dfrac{1}{n_{2}-\dfrac{h}{2}}\right ), \;h \text { is even, } \\ \prod\limits_{i=0}^{\frac{h-1}{2}}(n_{1}-i)(n_{2}-i), \;h \text { is odd. }\end{array}\right.$$
	\end{theorem}
	\begin{proof}
		Let $V(K_{n_{1},n_{2}})=N_{1}\cup N_{2}$ be a bipartition of $K_{n_{1},n_{2}}$, where $N_{1}=\left\{u_{1},\cdots,u_{n_{1}}\right\}$ and $N_{2}=\left\{v_{1},\cdots,v_{n_{2}}\right\}$.

		If $h$ is even, then the end-vertices of a path with length $h$ in $K_{n_{1},n_{2}}$ belong to the same part of its bipartition. Let
		$$u_{i_{0}}v_{i_{1}}u_{i_{2}}\cdots v_{i_{h-1}}u_{i_{h}}$$
be a path with length $h$ in $K_{n_{1},n_{2}}$, where $u_{i_{0}},u_{i_{2}},\cdots,u_{i_{h}}$ are $\dfrac{h}{2}+1$ vertices in $N_{1}$ and $v_{i_{1}},v_{i_{3}},\cdots, v_{i_{h-1}}$ are $\dfrac{h}{2}$ vertices in $N_{2}$. So the number of paths like $u_{i_{0}}v_{i_{1}}u_{i_{2}}\cdots v_{i_{h-1}}u_{i_{h}}$ is $$\dfrac{1}{2}\prod\limits_{i=0}^{\frac{h}{2}}(n_{1}-i)\prod\limits_{i=0}^{\frac{h}{2}-1}(n_{2}-i)=\dfrac{1}{2}\prod\limits_{i=0}^{\frac{h}{2}}(n_{1}-i)(n_{2}-i)\times\dfrac{1}{n_{2}-\dfrac{h}{2}}.$$

Similarly, the number of paths like ${v_{i_{0}}u_{i_{1}}v_{i_{2}}\cdots u_{j_{h-1}}v_{i_{h}}}$ is $$\dfrac{1}{2}\prod\limits_{i=0}^{\frac{h}{2}}(n_{2}-i)\prod\limits_{i=0}^{\frac{h}{2}-1}(n_{1}-i)=\dfrac{1}{2}\prod\limits_{i=0}^{\frac{h}{2}}(n_{2}-i)(n_{1}-i)\times\dfrac{1}{n_{1}-\dfrac{h}{2}}.$$

So, the number of paths with length $h$ in $K_{n_{1},n_{2}}$ is
		\begin{equation*}
			\begin{aligned} P_{h}(K_{n_{1},n_{2}})&=\dfrac{1}{2}\prod\limits_{i=0}^{\frac{h}{2}}(n_{1}-i)(n_{2}-i)\times\dfrac{1}{n_{2}-\dfrac{h}{2}}+\dfrac{1}{2}\prod\limits_{i=0}^{\frac{h}{2}}(n_{2}-i)(n_{1}-i)\times\dfrac{1}{n_{1}-\dfrac{h}{2}}\\
				&=\dfrac{1}{2}\prod\limits_{i=0}^{\frac{h}{2}}(n_{1}-i)(n_{2}-i)\times\left ( \dfrac{1}{n_{1}-\dfrac{h}{2}} +\dfrac{1}{n_{2}-\dfrac{h}{2}}\right ).
					\end{aligned}
			\end{equation*}

If $h$ is odd, then one of the end-vertices of the path belongs to $N_{1}$ and the other is contained in $N_{2}$. Let
		$$u_{i_{0}}v_{i_{1}}\cdots u_{i_{h-1}}v_{i_{h}}$$
be a path with length $h$ in $K_{n_{1},n_{2}}$, where $\{u_{i_{0}},u_{i_{2}},\cdots,u_{i_{h-1}}\}\subseteq N_1$ and $\{v_{i_{1}},v_{i_{3}},\cdots,v_{i_{h}}\}\subseteq N_2$. Then the number of paths like $u_{i_{0}}v_{i_{1}}\cdots u_{i_{h-1}}v_{i_{h}}$ is
$$\dfrac{1}{2}\prod\limits_{i=0}^{\frac{h-1}{2}}(n_{1}-i)\prod\limits_{i=0}^{\frac{h-1}{2}}(n_{2}-i)=\dfrac{1}{2}\prod\limits_{i=0}^{\frac{h-1}{2}}(n_{1}-i)(n_{2}-i).$$

Symmetrically, the number of paths like $v_{i_{0}}u_{i_{1}}\cdots v_{i_{h-1}}u_{i_{h}}$ is
$$\dfrac{1}{2}\prod\limits_{i=0}^{\frac{h-1}{2}}(n_{1}-i)\prod\limits_{i=0}^{\frac{h-1}{2}}(n_{2}-i)=\dfrac{1}{2}\prod\limits_{i=0}^{\frac{h-1}{2}}(n_{1}-i)(n_{2}-i).$$

So, the number of paths with length $h$ in $K_{n_1,n_2}$ is
		$$P_{h}(K_{n_{1},n_{2}})=\prod\limits_{i=0}^{\frac{h-1}{2}}(n_{1}-i)\prod\limits_{i=0}^{\frac{h-1}{2}}(n_{2}-i)=\prod\limits_{i=0}^{\frac{h-1}{2}}(n_{1}-i)(n_{2}-i).$$
\end{proof}


Note that a complete bipartite graph is determined by its bipartition, we can get the following result from the path sequence of a complete bipartite graph.

\begin{corollary}\label{c-5}
	Two complete bipartite graphs are isomorphic if and only if they have the same path sequence.
\end{corollary}
\begin{proof}
	If $K_{n_{1},n_{2}}$ and $K_{n^{\prime}_{1},n^{\prime}_{2}}$ are isomorphic, it is clearly that they have the same path sequence.
	
	If $K_{n_{1},n_{2}}$ and $K_{n^{\prime}_{1},n^{\prime}_{2}}$ have the same path sequence, i.e., $$P_{i}\left (K_{n_{1},n_{2}}\right )=P_{i}\left (K_{n^{\prime}_{1},n^{\prime}_{2}}\right ),i=0,\cdots,\rho.$$
	
From $$P_{0}\left (K_{n_{1},n_{2}}\right )=P_{0}\left (K_{n^{\prime}_{1},n^{\prime}_{2}}\right ),P_{1}\left (K_{n_{1},n_{2}}\right )=P_{1}\left (K_{n^{\prime}_{1},n^{\prime}_{2}}\right ), $$
we have $$n_{1}+n_2=n_{1}^{\prime}+n_2^{\prime}, n_1n_{2}=n_1^{\prime}n_{2}^{\prime}.$$
Without loss of generality, we may assume that $n_1\geq n_2$ and $n_1'\geq n_2'$. Since $f(x)=tx-x^2$ is monotonically decreasing for $x\geq\frac{t}{2}$, where $t=n_1+n_2$, we can get $n_1=n_1'$ and $n_2=n_2'$. It means that $K_{n_{1},n_{2}}$ and $K_{n^{\prime}_{1},n^{\prime}_{2}}$ are isomorphic.
\end{proof}

Specially, the path sequence of a star $S_n=K_{1,n-1}$ is
$$P(K_{1,n-1})=\left(n,n-1,\dfrac{(n-1)(n-2)}{2}\right).$$

For the path $P_{n}$ and the cycle $C_{n}$ with $n$ vertices, it is easy to get their path sequences. We have
$$P(P_{n})=\left(n,n-1,n-2,\cdots,1\right),$$
$$P(C_{n})=(\underbrace{n,n,\cdots,n}_{n}).$$

Now, we consider the number of paths of a starlike tree and their properties. Let $S_{n}$ be a starlike tree with $n$ vertices and root (i.e., center) $v_{0}$. We denote by $L_{l}$ the number of $l$-branches in $S_{n}$, where a $l$-branch is a path with length $l$ from the root $v_0$ to a pendant vertex, $l=1,\cdots,t$ and $t$ represents the length of a longest branch in $S_{n}$. Clearly, a starlike tree $S_n$ is completely determined its $L_1,L_2,\cdots,L_t$. In fact, Rada and Araujo [10] gave a method for calculating the number of paths in a starlike tree. This method is classified and calculated based on the sequence of degrees of vertices on the path.

 Let $\mathbb{N}=\{0,1,\cdots\}$ be the set of natural numbers. The function $\Psi_{h}: \mathcal{P}_{h}(G) \rightarrow \mathbb{N}^{h+1}$ is defined by $\Psi_{h}(\pi)=\left(d(v_{1}),d(v_{2}),\cdots,d(v_{h+1})\right)$, where $\pi=v_{1}v_{2}\cdots v_{h+1}$ is a path with length $h$ in $G$. So, a path in $\mathcal{P}_{h}(G)$ can be acted as a natural number sequence (i.e., degree sequence of the path) by the function $\Psi_{h}$. Thus, the number of paths is equivalent to the number of original images of degree sequences under $\Psi_{h}$.
	
From the proof of Theorem 2.3 in [10], we can get

\begin{lemma}\label{l-6} [10]
	The paths with length $h$ in $S_{n}$ can be counted in the following three types.\\
	\textbf{Type 1}
	The possible image under $\Psi_{h}$ of all paths in $\mathcal{P}_{h}(S_{n})$ which contain $v_{0}$ as an end-vertex is $X_{1}=(m, \underbrace{2,\cdots,2}_{h-1}, 1)$ or
	$X_{2}=(m, \underbrace{2,\cdots,2}_{h})$, and the number of paths in $\mathcal{P}_{h}(S_{n})$ whose images under $\Psi_{h}$ is $X_{1}$ or $X_{2}$ is
	$$ \left|\Psi_{h}^{-1}\left(X_{1}\right)\right|=L_{h},$$
	$$\left|\Psi_{h}^{-1}\left(X_{2}\right)\right|=L_{h+1}+\cdots+L_{t}=m-\sum_{i=1}^{h} L_{i}.$$
	\textbf{Type 2}
	The possible image under $\Psi_{h}$ of all paths in $ \mathcal{P}_{h}(S_{n})$ which do not contain $v_{0}$ is $Y_{1}=( \underbrace{2,\cdots,2}_{h+1})$ or $Y_{2}=(1, \underbrace{2,\cdots,2}_{h})$,  and the number of paths in $\mathcal{P}_{h}(S_{n})$ whose images under $\Psi_{h}$ is $Y_{1}$ or $Y_{2}$ is
	\begin{center}
		$$ \left|\Psi_{h}^{-1}\left(Y_{1}\right)\right|=(n-1)-\sum_{i=1}^{h}iL_{i}-(h+1)\left(m-\sum_{i=1}^{h}L_{i}\right),$$
		$$\left|\Psi_{h}^{-1}\left(Y_{2}\right)\right|=L_{h+1}+\cdots+L_{t}=m-\sum_{i=1}^{h}L_{i}.$$
	\end{center}
	\textbf{Type 3}
	The possible image under $\Psi_{h}$ of all paths in $ \mathcal{P}_{h}(S_{n})$ which contain $v_{0}$ but not as an end-vertex of the path is
	$$Z_{1}(a)=(1,\underbrace{2,\cdots,2}_{a},m,\underbrace{2,\cdots,2}_{h-1-a}), \quad 0\leq a\leq h-2$$ or
	$$Z_{2}(a)=(1, \underbrace{2, \cdots, 2}_{a}, m, \underbrace{2, \cdots, 2}_{h-2-a}, 1), \quad\left\{\begin{array}{l}0 \leq a \leq \dfrac{h}{2}-1, h \text { is even } \\ 0 \leq a \leq \dfrac{h-1}{2}-1, h \text { is odd }\end{array}\right.$$ or
	$$Z_{3}(a)=(\underbrace{2, \cdots, 2}_{a}, m, \underbrace{2, \cdots, 2}_{h-a}), \quad\left\{\begin{array}{l}1 \leq a \leq \dfrac{h}{2}, h \text { is even } \\ 1 \leq a \leq \dfrac{h-1}{2}, h \text { is odd}\end{array}\right.$$
	and the number of paths in $\mathcal{P}_{h}(S_{n})$ whose images under $\Psi_{h}$ is $Z_{1}$ or $Z_{2}$ or $Z_{3}$ is
	$$| \Psi_{h}^{-1}(Z_{1}(a)) |=\left\{\begin{array}{cc}\left\{\begin{array}{cc}L_{a+1}\left(m-\sum_{i=1}^{h-(a+1)}L_{i}\right) & \text { if } 0 \leqslant a \leqslant \frac{h}{2}-1 \\ L_{a+1}\left(m-1-\sum_{i=1}^{h-(a+1)} L_{i}\right) & \text { if } \frac{h}{2} \leqslant a \leqslant h-2\end{array}\right\} \text { if } h \text { is even, } \\ \left\{\begin{array}{cc}L_{a+1}\left(m-\sum_{i=1}^{h-(a+1)} L_{i}\right) & \text { if } 0 \leqslant a<\dfrac{h-1}{2} \\ L_{a+1}\left(m-1-\sum_{i=1}^{h-(a+1)} L_{i}\right) & \text { if } \dfrac{h-1}{2} \leqslant a \leqslant h-2\end{array}\right\} \text { if } h \text { is odd, }\end{array}\right.$$
	
	$$\left|\Psi_{h}^{-1}\left(\mathrm{Z}_{2}(a)\right)\right|=\left\{\begin{array}{c}\left\{\begin{array}{c}L_{a+1} L_{h-a-1} \quad \text { if } 0 \leqslant a<\dfrac{h}{2}-1 \\ \dfrac{1}{2}\left(L_{a+1}-1\right) L_{a+1} \quad \text { if } a=\dfrac{h}{2}-1\end{array}\right\} \text { if } h \text { is even, } \\ \left\{L_{a+1} L_{h-a-1} \quad \text { if } 0 \leqslant a \leqslant \dfrac{h-1}{2}-1\right\} \text { if } h \text { is odd, }\end{array}\right.$$
	$$\begin{array}{l}\left|\Psi_{h}^{-1}\left(\mathrm{Z}_{3}(a)\right)\right|=\left\{\begin{array}{l}\left\{\begin{array}{cc}\left[m-\sum_{i=1}^{h-a} L_{i}\right]\left[m-1-\sum_{i=1}^{a} L_{i}\right] \text { if } 1 \leqslant a \leqslant \frac{h}{2}-1 \\ \frac{1}{2}\left[m-\sum_{i=1}^{a} L_{i}\right]\left[m-1-\sum_{i=1}^{a} L_{i}\right] \quad \text { if } a=\frac{h}{2}\end{array}\right\}\text { if } h \text { is even, }\\ \left\{\left[m-\sum_{i=1}^{h-a} L_{i}\right]\left[m-1-\sum_{i=1}^{a} L_{i}\right] \quad \text { if } 1 \leqslant a \leqslant \frac{h-1}{2}\right\} \text { if } h \text { is odd. }\end{array}\right.\end{array}	$$
\end{lemma}	

By Lemma 5, we can show that the number of paths with length $h$ in a starlike tree is completely determined by its branches with length not more than $h-2$.

\begin{theorem}\label{t-7}
	If $S_{n}$ is a starlike tree with $n$ vertices and the degree $m$ of its root, then
	$$P_{h}(S_{n})=|\mathcal{P}_{h}(S_{n})|=\lambda(n,m,h,L_{1},\cdots,L_{h-3})+\mu(n,m,h,L_{1},\cdots,L_{h-3})L_{h-2}$$
	where $\lambda(n,m,h,L_{1},\cdots,L_{h-3})$ is a real number determined by the values of $n,m,h,L_{1},\cdots,L_{h-3}$ and $\mu(n,m,h,L_{1},\cdots,L_{h-3})=2-m, h>2$.
\end{theorem}
\begin{proof}
	From the definition of path sequence and Lemma 5, we can easily get
	$$P_{0}(S_{n})=n$$
	$$P_{1}(S_{n})=n-1$$
	$$P_{2}(S_{n})=\dfrac{1}{2}m^2-\dfrac{3}{2}m+n-1$$
	$$P_{3}(S_{n})=m^2+m+n-1+(2-m)L_{1}$$
	$$P_{4}(S_{n})=n-1-\dfrac{13}{2}m+\dfrac{3}{2}m^2+(\dfrac{7}{2}-2m)L_{1}+\dfrac{1}{2}L_{1}^{2}+(2-m)L_{2}$$
	If $h>3$ is odd,
	\begin{equation*}
		\begin{aligned}
			P_{h}(S_{n})=&L_{h}+m-\sum_{i=1}^{h}L_{i}+(n-1)-\sum_{i=1}^{h}iL_{i}-(h+1)\left(m-\sum_{i=1}^{h}L_{i}\right)+m-\sum_{i=1}^{h} L_{i}\\
			&\left.+L_{1}\left(m-\sum_{i=1}^{h-1}L_{i}\right)+L_{1}L_{h-1}+\sum_{a=\frac{h+1}{2}}^{h-2}L_{a+1}\left(m-\sum_{i=1}^{h-\left(a+1\right)}L_{i}-1\right )\right.\\
			&\left.+\sum_{a=1}^{\frac{h-3}{2}}\left[L_{a+1}\left(m-\sum_{i=1}^{h-(a+1)}L_{i}\right)+L_{a+1}L_{h-\left(a+1\right)}+\left(m-\sum_{i=1}^{h-a}L_{i}\right)\left(m-\sum_{i=1}^{a}L_{i}-1\right)\right]\right.\\
			&\left.+L_{\frac{h+1}{2}}\left(m-\sum_{i=1}^{\frac{h-1}{2}} L_{i}-1\right)+\left(m-\sum_{i=1}^{\frac{h+1}{2}} L_{i}\right)\left(m-\sum_{i=1}^{\frac{h-1}{2}} L_{i}-1\right)\right.\\
			=&\lambda\left(n,m,h,L_{1},L_{2},\cdots,L_{h-3}\right)+(2-m)L_{h-2}
		\end{aligned}
	\end{equation*}
	
	If $h>4$ is even, then
	\begin{equation*}
		\begin{aligned}
			P_{h}(S_{n})=&L_{h}+m-\sum_{i=1}^{h}L_{i}+(n-1)-\sum_{i=1}^{h} iL_{i}-(h+1)\left(m-\sum_{i=1}^{h}L_{i}\right)+m-\sum_{i=1}^{h} L_{i}\\
			&\left.+L_{1}\left(m-\sum_{i=1}^{h-1}L_{i}\right)+L_{1}L_{h-1}+\frac{1}{2}\left(m-\sum_{i=1}^{\frac{h}{2}}L_{i}\right)\left(m-\sum_{i=1}^{\frac{h}{2}}L_{i}-1\right)\right.\\
			&\left.+\sum_{a=1}^{\frac{h}{2}-2}\left[L_{a+1}\left(m-\sum_{i=1}^{h-(a+1)}L_{i}\right)+L_{a+1}L_{h-\left(a+1\right)}+\left(m-\sum_{i=1}^{h-a}L_{i}\right)\left(m-\sum_{i=1}^{a} L_{i}-1\right)\right]\right.\\
			&\left.+L_{\frac{h}{2}}\left(m-\sum_{i=1}^{\frac{h}{2}}L_{i}\right)+\frac{1}{2}L_{\frac{h}{2}}\left(L_{\frac{h}{2}}-1\right)+\left(m-\sum_{i=1}^{\frac{h}{2}+1}L_{i}\right)\left(m-\sum_{i=1}^{\frac{h}{2}-1}L_{i}-1\right)\right.\\
			=&\lambda\left(n,m,h,L_{1},L_{2},\cdots,L_{h-3}\right)+(2-m)L_{h-2}
		\end{aligned}
	\end{equation*}
	Thus, we can get for any $h>2$, $$P_{h}(S_{n})=|\mathcal{P}_{h}(S_{n})|=\lambda(n,m,h,L_{1},\cdots,L_{h-3})+\mu(n,m,L_{1},\cdots,L_{h-3})L_{h-2},$$ where $\lambda(n,m,h,L_{1},\cdots,L_{h-3})$ is a real number determined by the values of $n,m,h,L_{1},\cdots,L_{h-3}$ and $\mu(n,m,h,L_{1},\cdots,L_{h-3})=2-m$.
\end{proof}

As a consequence, we have	
\begin{theorem}\label{t-8}
	Two starlike trees $S_{n}$ and $S_{n^{\prime}}$ are isomorphic, if and only if, $S_{n}$ and $S_{n^{\prime}}$ have the same path sequence.
\end{theorem}
\begin{proof}
	If $S_{n}\cong S_{n^{\prime}}$, then it is clearly that $P_{h}(S_{n})=P_{h}(S_{n^{\prime}})$ for all $0\leq h\leq \rho(S_{n})=\rho(S_{n'})$ and $P(S_{n})=P(S_{n^{\prime}})$.
	
	On the other hand, let $P(S_{n})=P(S_{n^{\prime}})$, i.e., $P_{i}(S_{n})=P_{i}(S_{n^{\prime}})$ for all $0\leq i\leq \rho(S_{n})=\rho(S_{n'})$. Then
	$$n=P_{0}(S_{n})=P_{0}(S_{n^{\prime}})=n^{\prime}$$
	by the definition of the path sequence, and
	$$ m=m^{\prime}$$
	by $P_{2}(S_{n})=P_{2}(S_{n^{\prime}})$ in Theorem 6. Now, by $P_{3}(S_{n})=P_{3}(S_{n^{\prime}})$
	we can get
	$$\lambda(n,m)+(2-m)L_{1}=\lambda(n^{\prime},m^{\prime})+(2-m^{\prime})L_{1}^{\prime}.$$
	So,
	$$L_{1}=L_{1}^{\prime}.$$
	Next, applying Theorem 6 for $h=4$, we get
	$$\lambda(n,m,L_{1})+(2-m)L_{2}=\lambda(n^{\prime},m^{\prime},L_{1}^{\prime})+(2-m^{\prime})L_{2}^{\prime},$$
	and
	$$L_{2}=L_{2}^{\prime}.$$
	Note that the length $\rho(S_{n})$ of a longest path in $S_{n}$ is at least $t+1$, where $t$ is the length of a longest branch in $S_{n}$. If $\rho(S_{n})\geq t+2$, it is clear that we can continue this process by the Theorem 6 repeatedly to conclude that
	$$n=n^{\prime},m=m^{\prime},L_{i}=L_{i}^{\prime},i=1,2,\cdots,t.$$
	Hence, $S_{n}$ and $S_{n^{\prime}}$ are isomorphic.
	
	If $\rho(S_{n})=t+1$, $P(S_{n})=P(S_{n^{\prime}})$, it implies that $S_{n}$ and $S_{n^{\prime}}$ have only one branch with the longest length $t$ and the other branches with length $1$, i.e., $S_{n}$ and $S_{n^{\prime}}$ are the same broom. They are obviously isomorphic.
\end{proof}

Next, we consider kite graphs. The kite graph, denoted by $K_{n_{1}}^{n_{2}}$, is obtained by appending the complete graph $K_{n_{1}}$ to a pendant vertex $x$ of the path $P_{n_{2}}$ with $n_{2}$ vertices, where $|V(K_{n_{1}}^{n_{2}})|=n_{1}+n_{2}-1$ and $d_{K_{n_{1}}^{n_{2}}}(x)=n_{1}$. Actually, $K_{n_{1}}^{n_{2}}$ is the line graph of the starlike tree (a broom) $S_{n_{1}+n_{2}}$ with $L_{1}=n_1-1$, $L_{2}=\cdots=L_{n_{2}-1}=0$ and $L_{n_{2}}=1$. We can get the number of paths with length $h$ in $K_{n_{1}}^{n_{2}}$ for $0\leq h\leq\rho(K_{n_{1}}^{n_{2}})=n_{1}+n_{2}-2$.

 	\begin{theorem}\label{l-9}
 			For $0\leq h\leq n_{1}+n_{2}-2$, the number of paths with length $h$ in $K_{n_{1}}^{n_{2}}$
 $$P_{h}(K_{n_{1}}^{n_{2}})
 =\left\{\begin{array}{l}
 n_{1}+n_{2}-1, \; h =0\\
 \dfrac{n_{1}(n_{1}-1)}{2}+n_{2}-1, \;	h =1\\

 \left\{\begin{array}{l}
 \dfrac{1}{2}\prod\limits_{i=0}^{h}(n_{1}-i)+(n_{2}-h)+\sum\limits_{k=1}^{h-1}\prod\limits_{i=1}^{k}\left(n_{1}-i\right), \;2 \leq h \leq n_{2}-1 \\ \dfrac{1}{2}\prod\limits_{i=0}^{h}(n_{1}-i)+\sum\limits_{k=1}^{n_{2}-1}\prod\limits_{i=1}^{h-k}\left(n_{1}-i\right), \;n_{2}\leq h \leq n_{1}-1 \\ \sum\limits_{k=1}^{n_{2}-1}\prod\limits_{i=1}^{h-k}\left(n_{1}-i\right), \;n_{1}\leq	h \leq n_{1}+n_{2}-2 \\
 \end{array}\right\} \; n_2\leq n_1\\

\left\{\begin{array}{l}
 \dfrac{1}{2}\prod\limits_{i=0}^{h}(n_{1}-i)+(n_{2}-h)+\sum\limits_{k=1}^{h-1}\prod\limits_{i=1}^{k}\left(n_{1}-i\right), \;2 \leq h \leq n_{1}-1 \\ (n_2-h)+\sum\limits_{k=1}^{n_{1}-1}\prod\limits_{i=1}^{k}\left(n_{1}-i\right), \;n_{1}\leq h \leq n_{2}-1 \\
 \sum\limits_{k=1}^{n_{1}-1}\prod\limits_{i=1}^{k}\left(n_{1}-i\right), \;n_{2}\leq	h \leq n_{1}+n_{2}-2 \\
 \end{array}\right\} \; n_1< n_2\\

 \end{array}\right.$$
 			\end{theorem}

 \begin{proof}
 		First, we have
 		$$P_{0}(K_{n_{1}}^{n_{2}})=|V(K_{n_{1}}^{n_{2}})|=n_{1}+n_{2}-1,$$
 		$$P_{1}(K_{n_{1}}^{n_{2}})=|E(K_{n_{1}}^{n_{2}})|=\dfrac{n_{1}(n_{1}-1)}{2}+n_{2}-1.$$
 		
  		For $2\leq h\leq n_{1}+n_{2}-2$, all paths of length $h$ in $K_{n_{1}}^{n_{2}}$ can be divided into three types with the vertex $x$: the path contained entirely in $K_{n_{1}}$, the path contained entirely in $P_{n_{2}}$, and the path is connected by $x$ between a path of length $k$ in $K_{n_{1}}$ and a path of length $h-k$ in $P_{n_{2}}$, where $1\leq k\leq h-1$.
 		
Let $g(k)$ be the number of paths with length $k$ in $K_{n_{1}}$ which contain $x$ as an end-vertex, then $g(k)=\prod\limits_{i=1}^{k}\left(n_{1}-i\right)$. And the number of paths with length $k$ in $P_{n_{2}}$ which contain $x$ as an end-vertex is $1$.

 (I) $2\leq n_{2}\leq n_{1}$.

 If $2\leq h\leq n_{2}-1$, then we have
 		\begin{equation*}
 			\begin{aligned} P_{h}(K_{n_{1}}^{n_{2}})=P_{h}(K_{n_{1}})+P_{h}(P_{n_{2}})+\sum_{k=1}^{h-1}g(k)=\dfrac{1}{2}\prod_{i=0}^{h}(n_{1}-i)+(n_{2}-h)+\sum_{k=1}^{h-1}\prod_{i=1}^{k}\left(n_{1}-i\right).
 			\end{aligned}
 		\end{equation*}

If $n_{2}\leq h\leq n_{1}-1$, then
 		\begin{equation*}
 			\begin{aligned} P_{h}(K_{n_{1}}^{n_{2}})=P_{h}(K_{n_{1}})+\sum_{k=1}^{n_{2}-1}g(h-k)=\dfrac{1}{2}\prod_{i=0}^{h}(n_{1}-i)+\sum_{k=1}^{n_{2}-1}\prod_{i=1}^{h-k}\left(n_{1}-i\right).
 			\end{aligned}
 		\end{equation*}

 If $n_{1}\leq h\leq n_{1}+n_{2}-2$, then
 		\begin{equation*}
 			\begin{aligned} P_{h}(K_{n_{1}}^{n_{2}})=\sum_{k=1}^{n_{2}-1}g(h-k)=\sum_{k=-1}^{n_{2}-1}\prod_{i=1}^{k}\left(n_{1}-i\right).
 			\end{aligned}
 		\end{equation*}

 (II) $2\leq n_{1}< n_{2}$.

 If $2\leq h\leq n_{1}-1$, then
 		\begin{equation*}
 			\begin{aligned} P_{h}(K_{n_{1}}^{n_{2}})=P_{h}(K_{n_{1}})+P_{h}(P_{n_{2}})+\sum_{k=1}^{h-1}g(k)=\dfrac{1}{2}\prod_{i=0}^{h}(n_{1}-i)+(n_{2}-h)+\sum_{k=1}^{h-1}\prod_{i=1}^{k}\left(n_{1}-i\right).
 			\end{aligned}
 		\end{equation*}

If $n_{1}\leq h\leq n_{2}-1$, then
 		\begin{equation*}
 			\begin{aligned} P_{h}(K_{n_{1}}^{n_{2}})=P_{h}(P_{n_{2}})+\sum_{k=1}^{n_{1}-1}g(k)=(n_{2}-h)+\sum_{k=1}^{n_{1}-1}\prod_{i=1}^{k}\left(n_{1}-i\right).
 			\end{aligned}
 		\end{equation*}

 If $n_{2}\leq h\leq n_{1}+n_{2}-2$, then
 		\begin{equation*}
 			\begin{aligned} P_{h}(K_{n_{1}}^{n_{2}})=\sum_{k=1}^{n_{1}-1}g(k)=\sum_{k=1}^{n_{1}-1}\prod_{i=1}^{k}\left(n_{1}-i\right).
 			\end{aligned}
 		\end{equation*}

  		\end{proof}

\begin{corollary}\label{c-10}
	Two kite graphs $K_{n_{1}}^{n_{2}}$ and $K_{n_{1}^{\prime}}^{n_{2}^{\prime}}$ are isomorphic if and only if they have the same path sequence.
\end{corollary}		
 		
 \begin{proof}
 			If $K_{n_{1}}^{n_{2}}$ and $K_{n_{1}^{\prime}}^{n_{2}^{\prime}}$ are isomorphic, then it is clearly that they have the same path sequence.
 			
 			If $K_{n_{1}}^{n_{2}},K_{n_{1}^{\prime}}^{n_{2}^{\prime}}$ have the same path sequence, then
 $$P_{0}(K_{n_{1}}^{n_{2}})=P_{0}(K_{n_{1}^{\prime}}^{n_{2}^{\prime}}), \;\;\; P_{1}(K_{n_{1}}^{n_{2}})=P_{1}(K_{n_{1}^{\prime}}^{n_{2}^{\prime}}).$$
From Theorem 8, we have
 			$$n_{1}+n_{2}-1=n_{1}^{\prime}+n_{2}^{\prime}-1,\;\;\; \dfrac{n_{1}(n_{1}-1)}{2}+n_{2}-1=\dfrac{n_{1}^{\prime}(n_{1}^{\prime}-1)}{2}+n^{\prime}_{2}-1.$$
So, $\dfrac{n_{1}(n_{1}-1)}{2}-n_{1}=\dfrac{n_{1}^{\prime}(n_{1}^{\prime}-1)}{2}-n^{\prime}_{1}$. Note that $f(x)=\frac{1}{2}x^{2}-\frac{3}{2}x$ is monotonically increasing for $x\geq \frac{3}{2}$, we can get $n_{1}=n_{1}^{\prime}$. And $n_{2}=n_{2}^{\prime}$. It means that $K_{n_{1}}^{n_{2}},K_{n_{1}^{\prime}}^{n_{2}^{\prime}}$ are isomorphic.
\end{proof}

  Finally, we consider a lollipop graph, where we replace the complete graph $K_{n_{1}}$ in the kite graph $K_{n_{1}}^{n_{2}}$ with a cycle $C_{n_{1}}$ on $n_{1}$ vertices, i.e., the lollipop graph, denoted by $C_{n_{1}}^{n_{2}}$, is obtained by appending the cycle $C_{n_{1}}$ to a pendant vertex $x$ of the path $P_{n_{2}}$ with $n_{2}$ vertices, where $|V(C_{n_{1}}^{n_{2}})|=n_{1}+n_{2}-1$ and $d_{C_{n_{1}}^{n_{2}}}(x)=3$ and $\rho(C_{n_{1}}^{n_{2}})=n_{1}+n_{2}-2$.

 		\begin{theorem}\label{l-11}
 			For $0\leq h\leq n_{1}+n_{2}-2$, the number of paths with length $h$ in $C_{n_{1}}^{n_{2}}$ is
 $$P_{h}(C_{n_{1}}^{n_{2}})=\left\{\begin{array}{l}
 n_{1}+n_{2}-1, \;	h =0,1 \\

 \left\{\begin{array}{l}
 n_{1}+n_{2}+h-2, \;2 \leq	h \leq n_{2}-1 \\
 n_{1}+2n_{2}-2, \;n_{2}\leq h \leq n_{1}-1 \\
 2n_1+2n_{2}-2h-2, \;n_{1}\leq 	h \leq n_{1}+n_{2}-2 \\
 \end{array}\right\} \; n_2\leq n_1 \\

 \left\{\begin{array}{l}
 n_{1}+n_{2}+h-2, \;2 \leq	h \leq n_{1}-1 \\
 2n_{1}+n_{2}-h-2, \;n_{1}\leq h \leq n_{2}-1 \\
 2n_1+2n_{2}-2h-2, \;n_{2}\leq 	h \leq n_{1}+n_{2}-2 \\
 \end{array}\right\} \; n_1<n_2 \\

 \end{array}\right.$$
 		\end{theorem}
 		
 \begin{proof}
 			From the definition of path sequence, we have
 			$$P_{0}(C_{n_{1}}^{n_{2}})=|V(C_{n_{1}}^{n_{2}})|=n_{1}+n_{2}-1,$$
 			$$P_{1}(C_{n_{1}}^{n_{2}})=|E(C_{n_{1}}^{n_{2}})|=n_{1}+n_{2}-1.$$
 			
  			For $2\leq h\leq n_{1}+n_{2}-2$, all paths of length $h$ in $C_{n_{1}}^{n_{2}}$ can be divided into three types with the vertex $x$: the path contained entirely in $C_{n_{1}}$, the path contained entirely in $P_{n_{2}}$, and the path is connected by $x$ between a path of length $k$ in $C_{n_{1}}$ and a path of length $h-k$ in $P_{n_{2}}$, where $1\leq k\leq h-1$.

  The number of paths with length $k$ in $C_{n_{1}}$ which contain $x$ as an end-vertex is $2$. And the number of paths with length $k$ in $P_{n_{2}}$ which contain $x$ as an end-vertex is $1$.

 		(I) $2\leq n_{2}\leq n_{1}$.

 If $2\leq h\leq n_{2}-1$, then
 			\begin{equation*}
 				P_{h}(C_{n_{1}}^{n_{2}})=P_{h}(C_{n_{1}})+P_{h}(P_{n_{2}})+\sum_{k=1}^{h-1}2=n_{1}+(n_{2}-h)+2(h-1)=n_{1}+n_{2}+h-2.
 				\end{equation*}

 If $n_{2}\leq h\leq n_{1}-1$, then
 			\begin{equation*}
 				P_{h}(C_{n_{1}}^{n_{2}})=P_{h}(C_{n_{1}})+\sum_{k=1}^{n_{2}-1}2=n_{1}+2(n_{2}-1)=n_{1}+2n_{2}-2.
 				\end{equation*}

 			If $n_{1}\leq h\leq n_{1}+n_{2}-2$, then
 			\begin{equation*}
 				P_{h}(C_{n_{1}}^{n_{2}})=\sum_{k=h-n_2+1}^{n_{1}-1}2=2(n_1+n_{2}-h-1)=2n_1+2n_{2}-2h-2.
 				\end{equation*}

 (II) $n_1<n_2$.

 If $2\leq h\leq n_{1}-1$, then
 \begin{equation*}
 	 P_{h}(C_{n_{1}}^{n_{2}})=P_{h}(C_{n_{1}})+P_{h}(P_{n_{2}})+\sum_{k=1}^{h-1}2=n_{1}+(n_{2}-h)+2(h-1)=n_{1}+n_{2}+h-2.
 	\end{equation*}

 If $n_{1}\leq h\leq n_{2}-1$, then
 \begin{equation*}
 	 P_{h}(C_{n_{1}}^{n_{2}})=P_{h}(S_{n_{2}})+\sum_{k=1}^{n_{1}-1}2=(n_{2}-h)+2(n_{1}-1)=2n_{1}+n_{2}-h-2.
 	 \end{equation*}

 If $n_{2}\leq h\leq n_{1}+n_{2}-2$, then
 \begin{equation*}
 	 P_{h}(C_{n_{1}}^{n_{2}})=\sum_{k=h-n_1+1}^{n_{2}-1}2=2(n_1+n_{2}-h-1)=2n_1+2n_{2}-2h-2.
 	 \end{equation*}

         		\end{proof}

		\begin{corollary}\label{c-12}
 			Two lollipop graphs $C_{n_{1}}^{n_{2}}$ and $C_{n_{1}^{\prime}}^{n_{2}^{\prime}}$ are isomorphic if and only if they have the same path sequence.
 		\end{corollary}		
 		\begin{proof}
 			If $C_{n_{1}}^{n_{2}}$ and $C_{n_{1}^{\prime}}^{n_{2}^{\prime}}$ are isomorphic, then it is clearly that they have the same path sequence.
 			
 			If $C_{n_{1}}^{n_{2}},C_{n_{1}^{\prime}}^{n_{2}^{\prime}}$ have the same path sequence, and $n_{1},n_{2},n_{1}^{\prime},n_{2}^{\prime}\geq2$. From $P_{0}(C_{n_{1}}^{n_{2}})=P_{0}(C_{n_{1}^{\prime}}^{n_{2}^{\prime}})$, we have $n_{1}+n_{2}=n_{1}^{\prime}+n_{2}^{\prime}$. If $ n_{1}=n_{1}^{\prime}$, then $n_{2}=n_{2}^{\prime}$, and $C_{n_{1}}^{n_{2}}$ and $C_{n_{1}^{\prime}}^{n_{2}^{\prime}}$ are isomorphic. For $n_{1}\neq n_{1}^{\prime}$,
 			
(\uppercase\expandafter{\romannumeral1}) $ n_{2}< n_{1}$

(\romannumeral1) If $n_{2}^\prime< n_{1}^\prime$, we can assume that $n_1>n_1^\prime>n_2^\prime>n_2$. With $P_{h}(C_{n_{1}}^{n_{2}})=P_{h}(C_{n_{1}^{\prime}}^{n_{2}^{\prime}})$, we have
$$\left\{\begin{array}{ll}
n_{1}+n_{2}+h-2=n_{1}^{\prime}+n_{2}^{\prime}+h-2, 2\leq h\leq n_2-1,& (1a) \\
	n_{1}+2n_{2}-2=n_{1}^{\prime}+n_{2}^{\prime}+h-2, n_2\leq h\leq n_2^\prime-1, & (1b)  \\
	n_{1}+2n_{2}-2=n_{1}^{\prime}+2n_{2}^{\prime}-2, n_2^\prime\leq h\leq n_1^\prime-1, & (1c)\\
	n_{1}+2n_{2}-2=2n_1^{\prime}+2n_{2}^{\prime}-2h-2,n_1^\prime\leq h\leq n_1-1, & (1d)\\
	2n_1+2n_{2}-2h-2=2n_1^{\prime}+2n_{2}^{\prime}-2h-2,n_1\leq h\leq n_1+n_{2}-2. & (1e)
\end{array}\right.$$

From $n_{1}+n_{2}=n_{1}^{\prime}+n_{2}^{\prime}$ and (1c), we can get $n_2=n'_2$, a contradiction.

(\romannumeral2) If $n_{2}^\prime\geq n_{1}^\prime$, we can assume that $n_1>n_2^\prime\geq n_1^\prime>n_2$.
With $P_{h}(C_{n_{1}}^{n_{2}})=P_{h}(C_{n_{1}^{\prime}}^{n_{2}^{\prime}})$, we have
$$\left\{\begin{array}{ll}
	n_{1}+n_{2}+h-2=n_{1}^{\prime}+n_{2}^{\prime}+h-2, 2\leq h\leq n_2-1,& (2a) \\
	n_{1}+2n_{2}-2=n_{1}^{\prime}+n_{2}^{\prime}+h-2, n_2\leq h\leq n_1^\prime-1, & (2b)  \\
	n_{1}+2n_{2}-2=2n_1^{\prime}+n_{2}^{\prime}-h-2, n_1^\prime\leq h\leq n_2^\prime-1,& (2c) \\
	n_{1}+2n_{2}-2=2n_1^{\prime}+2n_{2}^{\prime}-2h-2,n_2^\prime\leq h\leq n_1-1, &(2d) \\
	2n_1+2n_{2}-2h-2=2n_1^{\prime}+2n_{2}^{\prime}-2h-2,n_1\leq h\leq n_1+n_{2}-2.&(2e) 
\end{array}\right.$$
From $n_{1}+n_{2}=n_{1}^{\prime}+n_{2}^{\prime}$ and (2b), it must be $n_2=n_1^\prime-1$, and $n_1=n_2^\prime+1$. Substituting them into (2c), we get $n_1^\prime\leq 1\leq n_2^\prime-1$, a contradiction.

(\uppercase\expandafter{\romannumeral2}) $ n_{2}> n_{1}$

(\romannumeral1) If $n_{2}^\prime> n_{1}^\prime$, we can assume that $n_1<n_1^\prime<n_2^\prime<n_2$.
With $P_{h}(C_{n_{1}}^{n_{2}})=P_{h}(C_{n_{1}^{\prime}}^{n_{2}^{\prime}})$, we have
$$\left\{\begin{array}{ll}
	n_{1}+n_{2}+h-2=n_{1}^{\prime}+n_{2}^{\prime}+h-2, 2\leq h\leq n_1-1,& (3a)  \\
	2n_{1}+n_{2}-h-2=n_{1}^{\prime}+n_{2}^{\prime}+h-2, n_1\leq h\leq n_1^\prime-1,& (3b)  \\
	2n_{1}+n_{2}-h-2=2n_1^{\prime}+n_{2}^{\prime}-h-2, n_1^\prime\leq h\leq n_2^\prime-1,& (3c) \\
	2n_{1}+n_{2}-h-2=2n_1^{\prime}+2n_{2}^{\prime}-2h-2,n_2^\prime\leq h\leq n_2-1,& (3d) \\
	2n_1+2n_{2}-2h-2=2n_1^{\prime}+2n_{2}^{\prime}-2h-2,n_2\leq h\leq n_1+n_{2}-2.& (3e) 
\end{array}\right.$$
From (3c) and $n_{1}+n_{2}=n_{1}^{\prime}+n_{2}^{\prime}$, we get $n_1=n'_1$, a contradiction.

(\romannumeral2) If $n_{2}^\prime\leq n_{1}^\prime$, we can assume that $n_1<n_2^\prime\leq n_1^\prime<n_2$.
With $P_{h}(C_{n_{1}}^{n_{2}})=P_{h}(C_{n_{1}^{\prime}}^{n_{2}^{\prime}})$, we have
$$\left\{\begin{array}{ll}
	n_{1}+n_{2}+h-2=n_{1}^{\prime}+n_{2}^{\prime}+h-2, 2\leq h\leq n_1-1, & (4a) \label{eqsystem16} \\
	2n_{1}+n_{2}-h-2=n_1^{\prime}+n_{2}^{\prime}+h-2, n_1\leq h\leq n_2^\prime-1, & (4b) \label{eqsystem17}  \\
	2n_{1}+n_{2}-h-2=n_{1}^{\prime}+2n_{2}^{\prime}-2, n_2^\prime\leq h\leq n_1^\prime-1, & (4c) \label{eqsystem18}\\
	2n_{1}+n_{2}-h-2=2n_1^{\prime}+2n_{2}^{\prime}-2h-2,n_1^\prime\leq h\leq n_2-1, & (4d) \label{eqsystem19}\\
	2n_1+2n_{2}-2h-2=2n_1^{\prime}+2n_{2}^{\prime}-2h-2,n_2\leq h\leq n_1+n_{2}-2. & (4e) \label{eqsystem20}
\end{array}\right.$$
From (4c) and $n_{1}+n_{2}=n_{1}^{\prime}+n_{2}^{\prime}$, we get $h=\frac{n_1}{2}$, it contradicts $n_1\leq h\leq n_2^\prime-1$.

(\uppercase\expandafter{\romannumeral3}) $ n_{1}= n_{2}$

(\romannumeral1) If $n_{2}^\prime> n_{1}^\prime$, we can assume that $n_1^\prime<n_1=n_2<n_2^\prime$.
With $P_{h}(C_{n_{1}}^{n_{2}})=P_{h}(C_{n_{1}^{\prime}}^{n_{2}^{\prime}})$, we have
$$\left\{\begin{array}{ll}
	n_{1}+n_{2}+h-2=n_{1}^{\prime}+n_{2}^{\prime}+h-2, 2\leq h\leq n_1^\prime-1, & (5a) \label{eqsystem21} \\
	n_{1}+n_{2}+h-2=2n_1^{\prime}+n_{2}^{\prime}-h-2, n_1^\prime\leq h\leq n_1-1, & (5b) \label{eqsystem22}  \\
	2n_1+2n_{2}-2h-2=2n_1^{\prime}+n_{2}^{\prime}-h-2, n_1\leq h\leq n_2^\prime-1, & (5c) \label{eqsystem23}\\
	2n_1+2n_{2}-2h-2=2n_1^{\prime}+2n_{2}^{\prime}-2h-2,n_2^\prime\leq h\leq n_2-1. & (5d) \label{eqsystem24}
\end{array}\right.$$
From (5b) and $n_{1}+n_{2}=n_{1}^{\prime}+n_{2}^{\prime}$, we get $h=\frac{n_1^{\prime}}{2}$, it contradicts $n_1^\prime\leq h\leq n_1-1$.

(\romannumeral2) If $n_{2}^\prime< n_{1}^\prime$, we can assume that $n_2^\prime<n_1= n_2<n_1^\prime$.
With $P_{h}(C_{n_{1}}^{n_{2}})=P_{h}(C_{n_{1}^{\prime}}^{n_{2}^{\prime}})$, we have
$$\left\{\begin{array}{ll}
	n_{1}+n_{2}+h-2=n_{1}^{\prime}+n_{2}^{\prime}+h-2, 2\leq h\leq n'_2-1, & (6a) \label{eqsystem25} \\
	n_{1}+n_{2}+h-2=n_{1}^{\prime}+2n_{2}^{\prime}-2, n'_2\leq h\leq n_1-1=n_2-1, & (6b) \label{eqsystem26}  \\
	2n_1+2n_{2}-2h-2=n_{1}^{\prime}+2n_{2}^{\prime}-2, n_1=n_2\leq h\leq n_1^\prime-1, & (6c) \label{eqsystem27}\\
	2n_1+2n_{2}-2h-2=2n_1^{\prime}+2n_{2}^{\prime}-2h-2, n_1^\prime\leq h\leq n'_1+n'_2-2. & (6d) \label{eqsystem28}
\end{array}\right.$$

From (6b) and $n_{1}+n_{2}=n_{1}^{\prime}+n_{2}^{\prime}$, we get $h=n'_2$, and then $n'_2=n_1-1=n_2-1$ since $n'_2\leq h\leq n_1-1=n_2-1$. So, $n_1^\prime=n_2^\prime+2=n_1+1=n_2+1$.

From (6c) and $n_{1}+n_{2}=n_{1}^{\prime}+n_{2}^{\prime}$, we get $h=\frac{n_1^\prime}{2}$, and then $n_1^\prime-1=n_1\leq \frac{n_1^\prime}{2}$ since $n_1=n_2\leq h\leq n_1^\prime-1$.
So, $n'_1\leq 2$ and $n_2^\prime\leq 0$, a contradiction.

(\romannumeral3) If $n_{2}^\prime= n_{1}^\prime$, then $n_1=n_2=n_1^\prime=n_2^\prime$ and $C_{n_{1}}^{n_{2}}$ and $C_{n_{1}^{\prime}}^{n_{2}^{\prime}}$ are isomorphic.

To sum up, $C_{n_{1}}^{n_{2}}$ and $C_{n_{1}^{\prime}}^{n_{2}^{\prime}}$ are isomorphic if $C_{n_{1}}^{n_{2}}$ and $C_{n_{1}^{\prime}}^{n_{2}^{\prime}}$ have the same path sequence.	
 			
 \end{proof}


\section{Another graph family determined by the path sequence}
 	
For which graph families $\mathcal{G}$, can each graph $G\in \mathcal{G}$ be distinguished by the path sequence, i.e., $G_1, G_2\in \mathcal{G}$ are isomorphic if and only they have the same path sequence? From the previous section, we know that this is true for the following family $\mathcal{G}$: (i) composed of complete graphs; (ii) composed of complete bipartite graphs; (iii) composed of starlike trees; (iv) composed of kite graphs; (v) composed of lollipop graphs. In this section, we will find more graph families $\mathcal{G}$ with this property.

We consider a generalized starlike tree and determine the path sequences of it.
	
A coalescence of two vertex disjoint graphs $G_{1}=(V_{1}, E_{1})$ and $G_{2}=(V_{2}, E_{2})$ with respect to vertex $u \in V_{1}$ and vertex $v \in V_{2}$, denoted by $G_{1}(u=v)G_{2}$, is obtained from the union of $G_{1}$ and $G_{2}$ by identifying $u$ and $v$.
	
	Specially, the coalescence $G_{n} = K_{n_{1}}(u_{0}=v_{0})S_{n_{2}}$ with respect to vertex $u_{0}\in K_{n_{1}}$ and the root $v_{0}$ of $S_{n_{2}}$ is called a generalized starlike tree with $n=n_1+n_2-1$ vertices, where $K_{n_{1}}$ and $S_{n_{2}}$ are the complete graph with $n_{1}$ vertices and a starlike tree with $n_{2}$ vertices, respectively. More specifically, renaming both $u_{0}$ and $v_{0}$ with $x_{0}$ in the coalescence $G_{n} = K_{n_{1}}(u_{0}=v_{0})S_{n_{2}}$, where $x_{0} \notin V(K_{n_{1}})\cup V(S_{n_{2}})$. The vertex set of $G_{n}$ is $\left(V(K_{n_{1}}) \backslash\{u_{0}\}\right) \cup\left(V(S_{n_{2}}) \backslash\{v_{0}\}\right) \cup\left\{x_{0}\right\}$ and $d(x_{0})=d_{G_{n}}(x_{0})=d_{K_{n_1}}(u_0)+d_{S_{n_{2}}}(v_{0})=n_{1}-1+m$, where $m=d_{S_{n_{2}}}(v_{0})$.
	
		The set of all generalized starlike trees on $n$ vertices with $d(x_{0})=r(r>2)$ is denoted by $\Omega_{n,r}$. In the following, we will show that every graph in $\Omega_{n,r}$ is determined by its path sequence.
	
	The path of length $h$ in $G_{n}$ is counted by the following three types: the path contains only vertices in $K_{n_{1}}$, the path contains only vertices in $S_{n_{2}}$, and the path is connected by $x_{0}$ between a path of length $k$ in $K_{n_{1}}$ and a path of length $h-k$ in $S_{n_{2}}$, where $1\leq k\leq h-1$.
	
	Let $g(k)$ be the number of paths with length $k$ in $K_{n_{1}}$ which contain $x_{0}$ as an end-vertex. Then $g(k)=\prod\limits_{i=1}^{k}\left(n_{1}-i\right)$.
	
	Let $f(h-k)$ be the number of paths with length $h-k$ in $S_{n_{2}}$ which contain $x_{0}$ as an end-vertex. Then $f(h-k)=L_{h-k}+m-\sum\limits_{i=1}^{h-k}L_{i}=m-\sum\limits_{i=1}^{h-k-1}L_{i}$ by Lemma 5.
	
	The path of length $h$ in $G_{n}$ can be obtained from the number of paths in $K_{n_{1}}$ and $S_{n_{2}}$, where $P_{h}(K_{n_{1}})=\frac{1}{2}\prod\limits_{i=0}^{h}\left(n_{1}-i\right)$ and $P_{h}(S_{n_{2}})=\lambda[n_2,m,h,L_{1},\cdots,L_{h-3}]+(2-m)L_{h-2}$ by Theorem 1 and Theorem 6.

	\begin{theorem}\label{l-14}
Let $G_{n} = K_{n_{1}}(u_{0}=v_{0})S_{n_{2}}$ be the coalescence of $K_{n_{1}}$ and $S_{n_{2}}$ as above. Then the number of paths with length $h$ in $G_{n}$ is

	{\footnotesize
	$P_{h}(G_{n})=$
	\begin{equation*}
	\left\{\begin{array}{l}
		n_1+n_2-1, \; h=0 \\ \\
		\frac{n_1(n_1-1)}{2}+n_2-1, \; h=1 \\ \\
		\frac{1}{2}m^2+(n_1-\frac{5}{2})m+\frac{1}{2}n_1^3-\frac{3}{2}n_1^2+n, \; h=2 \\ \\
		\left\{\begin{array}{l}
			\left.\begin{array}{l}
				\left.\begin{array}{l}
					P_h(K_{n_{1}})+P_h(S_{n_{2}})+\sum\limits_{k=1}^{h-1}\left [ g(k)f(h-k) \right ],\\ \qquad \qquad \qquad \qquad \qquad \qquad 3 \leq h \leq n_{1}-1 \\
					P_h(S_{n_{2}})+\sum\limits_{k=1}^{n_1-1}\left [ g(k)f(h-k) \right ],\\ \qquad \qquad \qquad \qquad \qquad \quad n_{1}\leq h \leq t_{1}+t_{2} \\
					\sum\limits_{k=1}^{n_1-1}\left [ g(k)f(h-k) \right ], \\ \qquad \qquad \qquad  t_{1}+t_{2}+1\leq h \leq t_{1}+n_{1}-1 \\
				\end{array}\right\} \; n_1-1\leq t_{1}<t_{1}+t_{2} \\
				
				\left.\begin{array}{l}
					P_h(K_{n_{1}})+P_h(S_{n_{2}})+\sum\limits_{k=1}^{h-1}\left [ f(k)g(h-k) \right ],\\ \qquad \qquad \qquad \qquad \qquad \qquad  \qquad 3 \leq h \leq t_{1} \\
					P_h(K_{n_{1}})+P_h(S_{n_{2}})+\sum\limits_{k=1}^{t_1}\left [ f(k)g(h-k) \right ],\\ \qquad \qquad \qquad \qquad \qquad t_{1}+1\leq h \leq n_{1}-1 \\
					P_h(S_{n_{2}})+\sum\limits_{k=1}^{t_1}\left [ f(k)g(h-k) \right ],\\ \qquad \qquad \qquad \qquad \qquad \quad n_{1}\leq 	h \leq t_{1}+t_{2} \\
					\sum\limits_{k=1}^{t_1}\left [ f(k)g(h-k) \right ], \\ \qquad \qquad \qquad t_{1}+t_{2}+1\leq h \leq t_{1}+n_{1}-1 \\
				\end{array}\right\} \; t_{2}\leq t_{1}<n_1-1\leq t_{1}+t_{2} \\
			\end{array}\right\} \; n_1-1\leq t_{1}+t_{2} \\
			
			\left.\begin{array}{l}
				P_h(K_{n_{1}})+P_h(S_{n_{2}})+\sum\limits_{k=1}^{h-1}\left [ f(k)g(h-k) \right ], \;3 \leq h \leq t_{1}\\
				P_h(K_{n_{1}})+P_h(S_{n_{2}})+\sum\limits_{k=1}^{t_1}\left [ f(k)g(h-k) \right ], \;t_{1}+1\leq h \leq t_{1}+t_{2} \\
				P_h(K_{n_{1}})+\sum\limits_{k=1}^{t_1}\left [ f(k)g(h-k) \right ], \;t_{1}+t_{2}+1\leq h \leq n_{1}-1 \\
				\sum\limits_{k=1}^{t_1}\left [ f(k)g(h-k) \right ], \;n_{1}\leq h \leq t_{1}+n_{1}-1 \\
			\end{array}\right\} \; t_{1}+t_{2}<n_1-1 \\
		\end{array}\right\} \; t_{2}<n_1-1 \\
		
		\left\{\begin{array}{l}
			P_h(K_{n_{1}})+P_h(S_{n_{2}})+\sum\limits_{k=1}^{h-1}\left [ g(k)f(h-k) \right ], \;3 \leq h \leq n_{1}-1 \\
			P_h(S_{n_{2}})+\sum\limits_{k=1}^{n_1-1}\left [ g(k)f(h-k) \right ], \;n_{1}\leq h \leq t_{1}+n_{1}-1 \\
			P_h(S_{n_{2}}), \;t_{1}+n_{1}\leq h \leq t_{1}+t_{2} \\
		\end{array}\right\} \; n_1-1\leq t_{2} \\
		
	\end{array}\right.
	\end{equation*}
}
where $t_1$ is the length of a longest branch in $S_{n_2}$ and $t_2$ is the length of a longest branch after deleting one longest branch in $S_{n_2}$.
		\end{theorem}

\begin{proof}
	From the definition of path sequence, we have
	$$P_{0}(G_{n})=|V(G_{n})|=n_1+n_2-1,$$
	$$P_{1}(G_{n})=|E(G_{n})|=\frac{n_1(n_1-1)}{2}+n_2-1,$$
	$$P_{2}(G_{n})=P_{2}(K_{n_1})+P_{2}(S_{n_2})+f(1)g(1)=\frac{1}{2}m^2+(n_1-\frac{5}{2})m+\frac{1}{2}n_1^3-\frac{3}{2}n_1^2+n.$$

	For $3 \leq h \leq \rho(G_{n})=\max\lbrace t_1+n_1-1,t_1+t_2\rbrace$, all paths of length $h$ in $G_{n}$ can be divided into three types with the vertex $x_0$: the path contained entirely in $K_{n_1}$, the path contained entirely in $S_{n_2}$, and the path is connected by $x_0$ between a path of length $k$ in $K_{n_1}$ and a path of length $h-k$ in $S_{n_2}$, where $1\leq k\leq h-1$.
	
	(\uppercase\expandafter{\romannumeral1}) If $\rho(G_{n})=t_1+n_1-1$, i.e., $n_1-1>t_2$.
	
	(i) $\rho(K_{n_1})\leq\rho(S_{n_2})$, i.e., $n_1-1\leq t_1+t_2$
	
	(a) $ n_1-1\leq t_{1}<t_{1}+t_{2}$
	
	If $3 \leq h \leq n_{1}-1$, then
		\begin{align}
	P_h(G_{n})=&P_h(K_{n_{1}})+P_h(S_{n_{2}})+\sum\limits_{k=1}^{h-1}\left[g(k)f(h-k \right] \nonumber \\
	=&\lambda(n_2,m,h,L_{1},\cdots,L_{h-3})+(2-m)L_{h-2}+\frac{1}{2}\prod\limits_{i=0}^{h}\left(n_{1}-i\right)\nonumber\\&+\sum\limits_{k=1}^{h-1}\left[\prod\limits_{i=1}^{k}\left(n_{1}-i\right)\left(m-\sum\limits_{i=1}^{h-k-1}L_{i}\right)\right] \nonumber\\
	=&(3-m-n_1)L_{h-2}+\lambda(n_2,m,h,L_{1},\cdots,L_{h-3})+\left(n_{1}-1\right)\left(m-\sum\limits_{i=1}^{h-3}L_{i}\right)\nonumber\\&+\sum\limits_{k=2}^{h-1}\left[\prod\limits_{i=1}^{k}\left(n_{1}-i\right)\left(m-\sum\limits_{i=1}^{h-k-1}L_{i}\right)\right]+\frac{1}{2}\prod\limits_{i=0}^{h}\left(n_{1}-i\right)\nonumber\\
	=&(3-m-n_1)L_{h-2}+\lambda_1(n_1,n_2,m,h,L_{1},\cdots,L_{h-3}).\tag{1-1}
		\end{align}
	
	If $n_{1}\leq h \leq t_{1}+t_{2}$, then
		\begin{align}
	P_h(G_{n})=&P_h(S_{n_{2}})+\sum\limits_{k=1}^{n_1-1}\left [ g(k)f(h-k) \right ]\nonumber \\
	=&\lambda(n_2,m,h,L_{1},\cdots,L_{h-3})+(2-m)L_{h-2}+\sum\limits_{k=1}^{n_1-1}\left[\prod\limits_{i=1}^{k}\left(n_{1}-i\right)\left(m-\sum\limits_{i=1}^{h-k-1}L_{i}\right)\right]\nonumber \\
	=&(3-m-n_1)L_{h-2}+\lambda(n_2,m,h,L_{1},\cdots,L_{h-3})+\left(n_{1}-1\right)\left(m-\sum\limits_{i=1}^{h-3}L_{i}\right)\nonumber \\&+\sum\limits_{k=2}^{n_1-1}\left[\prod\limits_{i=1}^{k}\left(n_{1}-i\right)\left(m-\sum\limits_{i=1}^{h-k-1}L_{i}\right)\right]\nonumber \\
	=&(3-m-n_1)L_{h-2}+\lambda_2(n_1,n_2,m,h,L_{1},\cdots,L_{h-3}).\tag{1-2}
        \end{align}

If $t_{1}+t_{2}+1\leq h \leq t_{1}+n_{1}-1$, then
		\begin{align}
   P_h(G_{n})=&\sum\limits_{k=1}^{n_1-1}\left [ g(k)f(h-k) \right]=\sum\limits_{k=1}^{n_1-1}\left[\prod\limits_{i=1}^{k}\left(n_{1}-i\right)\left(m-\sum\limits_{i=1}^{h-k-1}L_{i}\right)\right]\nonumber\\
              =&(n_1-1)\left(m-\sum\limits_{i=1}^{h-2}L_{i}\right)+\sum\limits_{k=2}^{n_1-1}\left[\prod\limits_{i=1}^{k}\left(n_{1}-i\right)\left(m-\sum\limits_{i=1}^{h-k-1}L_{i}\right)\right]\nonumber\\
              =&(1-n_1)L_{h-2}+\lambda_3(n_1,m,L_1,\cdots,L_{h-3}).\tag{1-3}
\end{align}

	(b) $t_{2}\leq t_{1}<n_1-1\leq t_{1}+t_{2}$
	
	If $3 \leq h \leq t_{1}$, then the number of paths with length $h$ is similar to (1-1)
	\begin{align}
    P_h(G_{n})=&P_h(K_{n_{1}})+P_h(S_{n_{2}})+\sum\limits_{k=1}^{h-1}\left [ f(k)g(h-k) \right]\nonumber\\=&(3-m-n_1)L_{h-2}+\lambda_1(n_1,n_2,m,h,L_{1},\cdots,L_{h-3}).\tag{2-1}
\end{align}

	If $t_{1}+1\leq h \leq n_{1}-1 $, then
	\begin{align}
	P_h(G_{n})=&P_h(K_{n_{1}})+P_h(S_{n_{2}})+\sum\limits_{k=1}^{t_1}\left [ f(k)g(h-k) \right]\nonumber\\
=&(2-m)L_{h-2}+\lambda(n_2,m,h,L_{1},\cdots,L_{h-3})+\frac{1}{2}\prod\limits_{i=0}^{h}\left(n_{1}-i\right)\nonumber\\
&+\sum\limits_{k=1}^{t_1}\left[\prod\limits_{i=1}^{h-k}\left(n_{1}-i\right)\left(m-\sum\limits_{i=1}^{k-1}L_{i}\right)\right] \nonumber\\
=&\mu_2(n_1,m)L_{h-2}+\lambda_4(n_2,m,h,L_{1},\cdots,L_{h-3}).\tag{2-2}
\end{align}

	If $n_{1}\leq h \leq t_{1}+t_{2}$, then
	\begin{align}
	P_h(G_{n})=&P_h(S_{n_{2}})+\sum\limits_{k=1}^{t_1}\left [ f(k)g(h-k)\right]\nonumber\\=&(2-m)L_{h-2}+\lambda(n_2,m,h,L_{1},\cdots,L_{h-3})+\sum\limits_{k=1}^{t_1}\left[\prod\limits_{i=1}^{h-k}\left(n_{1}-i\right)\left(m-\sum\limits_{i=1}^{k-1}L_{i}\right)\right] \nonumber\\=&\mu_2(n_1,m)L_{h-2}+\lambda_5(n_1,n_2,m,h,L_{1},\cdots,L_{h-3}).\tag{2-3}
\end{align}

	If $t_{1}+t_{2}+1\leq h \leq t_{1}+n_{1}-1$, then
	\begin{align}
	P_h(G_{n})=&\sum\limits_{k=1}^{t_1}\left [ f(k)g(h-k) \right ]=\sum\limits_{k=1}^{t_1}\left[\prod\limits_{i=1}^{h-k}\left(n_{1}-i\right)\left(m-\sum\limits_{i=1}^{k-1}L_{i}\right)\right]\nonumber\\
=&\mu_3(n_1,m)L_{h-2}+\lambda_6(n_2,m,h,L_{1},\cdots,L_{h-3}),\tag{2-4}
\end{align}
where $\mu_3(n_1,m)=0$.

	(ii) $ \rho(K_{n_1})>\rho(S_{n_2})$
	
	If $3 \leq h \leq t_{1} $, then the number of paths with length $h$ is similar to (1-1)
	\begin{align}
	P_h(G_{n})=&P_h(K_{n_{1}})+P_h(S_{n_{2}})+\sum\limits_{k=1}^{h-1}\left [ f(k)g(h-k) \right ]\nonumber\\=&(3-m-n_1)L_{h-2}+\lambda_1(n_1,n_2,m,h,L_{1},\cdots,L_{h-3}).\tag{3-1}
\end{align}

	If $t_{1}+1\leq h \leq t_{1}+t_{2}$, then the number of paths with length $h$ is similar to (2-2)
	\begin{align}
	P_h(G_{n})=&P_h(K_{n_{1}})+P_h(S_{n_{2}})+\sum\limits_{k=1}^{t_1}\left [ f(k)g(h-k) \right ]\nonumber\\
=&(2-m)L_{h-2}+\lambda(n_2,m,h,L_{1},\cdots,L_{h-3})+\frac{1}{2}\prod\limits_{i=0}^{h}\left(n_{1}-i\right)\nonumber\\
&+\sum\limits_{k=1}^{t_1}\left[\prod\limits_{i=1}^{h-k}\left(n_{1}-i\right)\left(m-\sum\limits_{i=1}^{k-1}L_{i}\right)\right] \nonumber\\
=&\mu_2(n_1,m)L_{h-2}+\lambda_4(n_2,m,h,L_{1},\cdots,L_{h-3}).\tag{3-2}
\end{align}
	
	If $t_{1}+t_{2}+1\leq h \leq n_{1}-1$, then
	\begin{align}
    P_h(G_{n})=&P_h(K_{n_{1}})+\sum\limits_{k=1}^{t_1}\left [ f(k)g(h-k) \right ]\nonumber\\
    =&\frac{1}{2}\prod\limits_{i=0}^{h}\left(n_{1}-i\right)+\sum\limits_{k=1}^{t_1}\left[\prod\limits_{i=1}^{h-k}\left(n_{1}-i\right)\left(m-\sum\limits_{i=1}^{k-1}L_{i}\right)\right]\nonumber\\
    =&\mu_3(n_1,m)L_{h-2}+\lambda_7(n_2,m,h,L_{1},\cdots,L_{h-3}). \tag{3-3}
\end{align}

	If $n_{1}\leq h \leq t_{1}+n_{1}-1$, then the number of paths with length $h$ is similar to (2-4)
	\begin{align}
	P_h(G_{n})=&\sum\limits_{k=1}^{t_1}\left [ f(k)g(h-k) \right ]=\sum\limits_{k=1}^{t_1}\left[\prod\limits_{i=1}^{h-k}\left(n_{1}-i\right)\left(m-\sum\limits_{i=1}^{k-1}L_{i}\right)\right]\nonumber\\
=&\mu_3(n_1,m)L_{h-2}+\lambda_6(n_2,m,h,L_{1},\cdots,L_{h-3}).\tag{3-4}
\end{align}

(\uppercase\expandafter{\romannumeral2}) If $\rho(G_{n})=t_1+t_2$, i.e., $n_1-1\leq t_{2}$, then $\rho(K_{n_1})\leq \rho(S_{n_2})$.
	
	If $3 \leq h \leq n_{1}-1$, then the number of paths with length $h$ is similar to (1-1)
	\begin{align}
		P_h(G_{n})=&P_h(K_{n_{1}})+P_h(S_{n_{2}})+\sum\limits_{k=1}^{h-1}\left [ g(k)f(h-k) \right ]\nonumber\\=&(3-m-n_1)L_{h-2}+\lambda_1(n_1,n_2,m,h,L_{1},\cdots,L_{h-3}).\tag{4-1}
	\end{align}
	
	If $n_{1}\leq h \leq t_{1}+n_{1}-1$, then the number of paths with length $h$ is similar to (1-2)
		\begin{align}
		P_h(G_{n})=&P_h(S_{n_{2}})+\sum\limits_{k=1}^{n_1-1}\left [ g(k)f(h-k)\right]\nonumber\\=&(3-m-n_1)L_{h-2}+\lambda_2(n_1,n_2,m,h,L_{1},\cdots,L_{h-3}).\tag{4-2}
	\end{align}
	
	If $t_{1}+n_{1}\leq h \leq t_{1}+t_{2}$, then
	\begin{align}
	P_h(G_{n})=P_h(S_{n_{2}})=(2-m)L_{h-2}+\lambda(n,m,h,L_{1},\cdots,L_{h-3}).\tag{4-3}
\end{align}
	\end{proof}

From the proof of Theorem 12, we can get the following recurrent relation on the number of paths with length $h$ for the generalized starlike trees, which shows that the number of paths with length $h$ in the generalized starlike tree $G_n= K_{n_{1}}(u_{0}=v_{0})S_{n_{2}}$ is also completely by the branches with length not more than $h-2$ in $S_{n_2}$.
	
	\begin{corollary}\label{l-15}
		If $G_{n}= K_{n_{1}}(u_{0}=v_{0})S_{n_{2}}$, then
		$$P_{h}(G_{n})=\overline{\lambda}(n_1,n_2,m,h,L_{1},\cdots,L_{h-3})+\overline{\mu}(n_1,m)L_{h-2}$$
		where $L_i$ is the number of branches with length $i$ in $S_{n_2}$, $\overline{\lambda}(n_1,n_2,m,h,L_{1},\cdots,L_{h-3})$ is a real number determined by the values of $n_1, n_2, m$, $h, L_{1}, \cdots, L_{h-3}$, $\overline{\mu}(n_1,m)$ is a real number determined by the values of $n_1, m$; $\overline{\mu}(n_1,m)=3-m-n_1$ or $\overline{\mu}(n_1,m)=2-m$ for $2<h\leq t_1+t_2$.
	\end{corollary}

Note that a starlike tree $S_n$ is completely determined by its sequence $L_1,L_2,\cdots,L_t$ of branches, where $L_i$ is the number of branches with length $i$ in $S_n$ and $t$ is the length of a longest branch. From the definition of the coalescence, it is easy to know the following fact:

Two generalized starlike trees $G_{n}= K_{n_{1}}(u_{0}=v_{0})S_{n_{2}}$ and $G_{n'}= K_{n'_{1}}(u'_{0}=v'_{0})S_{n'_{2}}$ are isomorphic, if and only, $n_1=n'_1$ and $S_{n_2}$ and $S_{n'_2}$ have the same sequence of branches (it implies $m=m'$).

Note that $m=\sum\limits_{i=1}^{t}L_i$, $m'=\sum\limits_{i=1}^{t'}L'_i$, $n_2=1+\sum\limits_{i=1}^{t}iL_i$ and $n'_2=1+\sum\limits_{i=1}^{t'}iL'_i$, if $m=m'$, $n_2=n'_2$ and $L_1=L'_1, L_2=L'_2, \cdots, L_{t-1}=L'_{t-1}$, without loss of generality, assuming $t\leq t'$, we have
$$L_t=L'_t+L'_{t+1}+\cdots,\,\,\,\, \,\,\, tL_t=tL'_t+(t+1)L'_{t+1}+\cdots.$$
Then $L'_{t+1}=L'_{t+2}=\cdots=0$ and $L_t=L'_t$. And $S_{n_2}$ and $S_{n'_2}$ have the same sequence of branches.

	\begin{theorem}\label{t-16}
		Let $G_{n}= K_{n_{1}}(u_{0}=v_{0})S_{n_{2}}$ and $G_{n^{\prime}}= K_{n'_{1}}(u'_{0}=v'_{0})S_{n'_{2}}$ be two generalized starlike trees. Then $G_{n}$ and $G_{n}^{\prime}$ are isomorphic, if and only if, they have the same path sequence.
	\end{theorem}
	\begin{proof}
		Let $m$ and $m'$ be the degrees of roots in $S_{n_2}$ and $S_{n'_2}$, $t_1,t'_1$ the length of a longest branch in $S_{n_2}$,$S_{n'_2}$ and $t_2,t'_2$ be the length of the longest branch after deleting a longest branch in $S_{n_2}$, $S_{n'_2}$, respectively.
		
		If $G_{n}\cong G_{n^{\prime}}$, then they have clearly the same path sequence.
		
		On the other hand, we assume they have the same path sequence, i.e., $P_{h}(G_{n})=P_{h}(G_{n^{\prime}})$ for all $0\leq h\leq \rho(G_{n})=\rho(G_{n'})$. In order to show that $G_{n}$ and $G_{n}^{\prime}$ are isomorphic, we only need to show that $n_1=n'_1$, $n_2=n'_2$, $m=m'$ and $L_1=L'_1, L_2=L'_2, \cdots, L_{t_1-1}=L'_{t_1-1}$ by the fact above.

Using Theorem 12, we have
		$$n_{1}+n_{2}-1=n_{1}^{\prime}+n_{2}^{\prime}-1,\;\;\; \dfrac{n_{1}(n_{1}-1)}{2}+n_{2}-1=\dfrac{n_{1}^{\prime}(n_{1}^{\prime}-1)}{2}+n^{\prime}_{2}-1$$
by $P_{0}(G_{n})=P_{0}(G_{n^{\prime}})$ and $P_{1}(G_{n})=P_{1}(G_{n^{\prime}})$, respectively. So, $\frac{1}{2}n_1^{2}-\frac{3}{2}n_1=\frac{1}{2}n'^{2}_1-\frac{3}{2}n'_1$.

	Note that $f(x)=\frac{1}{2}x^{2}-\frac{3}{2}x$ is monotonically increasing for $x\geq \frac{3}{2}$, we can get $n_{1}=n_{1}^{\prime}$ and $n_{2}=n_{2}^{\prime}$. Moreover, $m=m^{\prime}$ by $P_{2}(G_{n})=P_{2}(G_{n^{\prime}})$.

From Corollary 13, we have
$$\overline{\lambda}(n_1,n_2,m,h,L_{1},\cdots,L_{h-3})+\overline{\mu}(n_1,m)L_{h-2}=\overline{\lambda}(n'_1,n'_2,m',h,L'_{1},\cdots,L'_{h-3})+\overline{\mu}(n'_1,m')L'_{h-2}$$
by using $P_{h}(G_{n})=P_{h}(G_{n'})$ for $h\geq 3$.

Recursively, we can obtain
$$L_1=L'_1, L_2=L'_2, \cdots, L_{t_1+t_2}=L'_{t_1+t_2}.$$
As a consequence, $G_{n}$ and $G_{n}^{\prime}$ are isomorphic.
			
	\end{proof}

\section {Conclusions}
	
	In this paper, we introduce the concept of the path sequence in a graph, and consider whether a graph can be determined by its path sequence. First, we obtain the number of paths with length $h$ in some graphs and prove that the graph in these graph families can be determined by the path sequence. Further more, we construct a new graph family, composed of the generalized starlike trees $G_{n}=K_{n_{1}}(u_{0}=v_{0})S_{n_{2}}$. We also get its path sequence and prove that $G_{n}$ can also be determined by its path sequence. We expect to find other graph families that can be determined by the path sequence. The following topics are worth further exploring:

Find more graph families where each graph can be distinguished by a path sequence. Or, construct some simple and connected graphs such that they have the same sequence.

\section*{Acknowledgements}
This work is supported by the National Natural Science Foundation of China (No.12201634) and the Hunan Provincial Natural Science Foundation of China (2023JJ30070,2020JJ4423).

\end{document}